\documentclass[twosided,reqno]{amsart}

\usepackage{amsmath}
\usepackage{amssymb}
\usepackage{amsthm}

\theoremstyle{theorem}
\newtheorem{theorem}{\scshape Theorem }[section]

\newtheorem{corollary}[theorem]{\scshape Corollary}

\theoremstyle{definition}

\newtheorem*{remark*}{\scshape Remark}
\newcommand{\ma}{\mathbb}
\newcommand{\be}{\begin{equation}}
\newcommand{\ee}{\end{equation}}
\newcommand{\ben}{\begin{equation*}}
\newcommand{\een}{\end{equation*}}
\newcommand{\bc}{\begin{corollary}}
\newcommand{\ec}{\end{corollary}}
\newcommand{\fa}{\frac}
\newcommand{\la}{\label}
\newcommand{\U}{\sum_{n=0}^{\infty}}
\newcommand{\bt}{\begin{theorem}}
\newcommand{\et}{\end{theorem}}
\newcommand{\bi}{\binom}
\newcommand{\lp}{\left(}
\newcommand{\rp}{\right)}
\newcommand{\tn}{\frac{t^n}{n!}}
\newcommand{\iq}{\lim_{q\to 1}}
\newcommand{\N}{\mathbb{N}}
\newcommand{\zp}{\mathbb{Z}_p}
\newcommand{\cp}{\mathbb{C}_p}
\newcommand{\eq}{E_{n,q}^{(h,r)}(x)}
\newcommand{\er}{E_{n,q}^{(r)}(x)}
\newcommand{\en}{E_{n}^{(r)}(x)}
\newcommand{\bq}{(b;q)_n}
\newcommand{\BQ}{(1-b)(1-bq)\cdots (1-bq^{n-1})}
\newcommand{\ms}{m_1, \cdots ,m_r=0}
\newcommand{\ys}{y_1, \cdots ,y_r=0}
\newcommand{\YS}{y_1+ \cdots +y_r}
\newcommand{\js}{j_1, \cdots ,j_r=0}
\newcommand{\JS}{j_1+ \cdots +j_r}
\newcommand{\is}{i_1, \cdots, i_r=0}
\newcommand{\MS}{m_1+ \cdots +m_r}
\newcommand{\ty}{\infty}
\newcommand{\I}{\int_{\mathbb{Z}_p} \cdots \int_{\mathbb{Z}_p}}
\newcommand{\dm}{d\mu_{-1}(y_1) \cdots d\mu_{-1}(y_r)}
\newcommand{\LN}{\lim_{N\to \infty}}
\newcommand{\iz}{\int_{\mathbb{Z}_p}}
\newcommand{\lr}{\sum_{l=1}^{r}}

\numberwithin{equation}{section}

\begin{document}

\title{Identities of symmetry for expansions of $q$-Euler polynomials}

\author{Dae San Kim}
\address{Department of Mathematics, Sogang University, Seoul 121-742, Republic of Korea.}
\email{dskim@sogang.ac.kr}

\author{Tae Gyun Kim}
\address{Department of Mathematics, Kwangwoon University, Seoul 139-701, Republic of Korea}
\email{tgkim2013@hotmail.com}

\maketitle

\begin{abstract}
Recently, Kim considered expansions of $q$-Euler polynomials which are given by
\ben
\U\eq\tn=2^r \sum_{\ms}^{\ty}q^{\lr(h-l)m_l}(-1)^{\lr m_l}e^{[x+\lr m_l]_q t}.
\een
In this paper, we investigate some symmetric properties of the multivariate $p$-adic fermionic integrals on $\zp$ and derive various identities concerning the expansions of $q$-Euler polynomials from the symmetric properties of the multivariate $p$-adic fermionic integrals on $\zp$
\end{abstract}

\section{Introduction}

Let $p$ be a fixed odd prime number. Throughout this paper $\zp, \ma{Q}_p, \cp$ will, respectively, denote the ring of $p$-adic rational integers, the field of $p$-adic rational numbers and the completion of algebraic closure of $\ma{Q}_p$. Let $\nu_p$ be the normalized exponential valuation of $\cp$ with $|p|_p=p^{-\nu_p(p)}=p^{-1}$. When one talks about a $q$-extension, $q$ is variously considered as an indeterminate, a complex number $q \in \ma{C}$ or a $p$-adic number $q\in \cp$. If $q \in \ma{C}$, one usually assumes $|q|<1$;if $q\in \cp$, one usually assumes $|1-q|_p<1$. The $q$-number of $x$ is defined by $[x]_q=\fa{1-q^x}{1-q}$. Note that $\iq[x]_q=x$. The $q$-factorial is defined as $[n]_q!=[n]_q[n-1]_q \cdots [2]_q[1]_q$. As is well known, the $q$-binomial formulae are given by

\be\la{1}
\begin{split}
\bq&=\BQ\\
&=\sum_{i=0}^{n}\bi{n}{i}_q q^{\bi{i}{2}}(-1)^i b^i,~~\textrm{(see $[1,5]$)},\\
\end{split}
\ee
where $\bi{n}{i}_q=\fa{[n]_q!}{[i]_q![n-i]_q!}=\fa{[n]_q\cdots [n-i+1]_q}{[i]_q!}$, and

\be\la{2}
\fa{1}{\bq}=\fa{1}{\BQ}=\sum_{i=0}^{\ty}\bi{n+i-1}{i}_q b^i.
\ee

Let $C(\zp)$ be the space of continuous functions on $\zp$. For $f\in C(\zp)$, the $p$-adic fermionic integral on $\zp$ is defined by Kim as follows:

\be\la{3}
I_{-1}(f)=\iz f(x)d\mu_{-1}(x)=\LN \sum_{x=0}^{p^N-1}f(x)(-1)^x, ~~\textrm{(see \cite{06})}.
\ee

From (\ref{3}), we note that

\be\la{4}
I_{-1}(f_n)+(-1)^{n-1}I_{-1}(f)=2\sum_{l=0}^{n-1}(-1)^{n-l-1}f(l),
\ee
where $f_n(x)=f(x+n)$ and $n\in \N$.\\
As is well known, the higher-order Euler polynomials are defined by the generating function to be

\be\la{5}
\lp\fa{2}{e^t+1}\rp^r e^{xt}=\U \en\tn,~(r\in \N), ~~\textrm{(see [1-14])}.
\ee

The $q$-extension of (\ref{5}) is given by

\be\la{6}
\begin{split}
\U\er\tn&=2^r\sum_{\ms}^{\ty}(-1)^{\MS}e^{[\MS+x]_qt}\\
&=2^r\sum_{m=0}^{\ty}\bi{m+r-1}{m}(-1)^me^{[m+x]_q t}, ~~\textrm{(see [5,6])}.
\end{split}
\ee

Recently, Kim considered expansions of $\er$ which are given by

\be\la{7}
\begin{split}
\U\eq\tn&=2^r\sum_{\ms}^{\ty}q^{\lr (h-l)m_l}(-1)^{\lr m_l}e^{[x+\MS]_q t}\\
&=2^r\sum_{m=0}^{\ty}\bi{m+r-1}{m}_q q^{(h-r)m}(-1)^m e^{[m+x]_q t},
\end{split}
\ee
where $h\in\ma{Z}$ and $r\in\N$ (see [5,6,7]).\\
In this paper, we investigate some symmetric properties of the multivariate $p$-adic fermionic integrals on $\zp$ and derive various identities concerning the expansions of $q$-Euler polynomials from the symmetric properties of the multivariate $p$-adic fermionic interals on $\zp$.

\section{Identities of Symmetry for expansions of $q$-Euler polynomials}

For $h \in \ma{Z}$ and $r\in \N$, from (\ref{3}), we have

\be\la{8}
\begin{split}
&\I q^{\lr(h-l)y_l}e^{[x+\YS]_q t}\dm\\
&=2^r \sum_{\ms}^{\ty}q^{\lr (h-l)m_l}(-1)^{\lr m_l}e^{[x+\MS]_q t}\\
&=2^r\sum_{m=0}^{\ty}\bi{m+r-1}{m}_q q^{(h-r)m}(-1)^m e^{[m+x]_q t}.
\end{split}
\ee

Thus, by (\ref{7}) and (\ref{8}), we get

\be\la{9}
\begin{split}
\eq&=\I q^{\lr(h-l)y_l}[x+\YS]_q^n\dm\\
&=2^r\sum_{m=0}^{\ty}\bi{m+r-1}{m}_q (-q^{h-r})^m[x+m]_q^n\\
&=\fa{2^r}{(1-q)^n}\sum_{l=0}^{n}\fa{\bi{n}{l}(-q^x)^l}{(-q^{h-r+l};q)_r}\\
&=2^r \sum_{\ms}^{\ty}q^{\lr(h-l)m_l}(-1)^{\lr m_l}[\MS+x]_q^n,
\end{split}
\ee
where $n \geq 0$.\\

Let $w_1, w_2 \in \N$ with $w_1 \equiv 1\pmod{2}$ and $w_2 \equiv 1\pmod{2}$. Then we see that

\be\la{10}
\begin{split}
&\I q^{w_1 \lr(h-l)y_l} e^{[w_1]_q[w_2x+\fa{w_2}{w_1}\lr j_l+\lr y_l]_{q^{w_1}} t}\dm\\
&=\I q^{w_1 \lr(h-l)y_l} e^{[w_1 w_2 x+w_2\lr j_l +w_1 \lr y_l]_q t}\dm\\
&=\LN \sum_{\ys}^{w_2p^{N}-1} q^{w_1 \lr(h-l)y_l} e^{[w_1 w_2 x+w_2\lr j_l +w_1 \lr y_l]_q t}(-1)^{\lr y_l}\\
&\hspace{3cm}\times\dm\\
&=\LN \sum_{\is}^{w_2-1}\sum_{\ys}^{p^{N}-1} q^{w_1 \lr(h-l)(i_l+w_2y_l)}\\
&\hspace{3cm}\times e^{[w_1 w_2 x+w_2\lr j_l +w_1 \lr (i_l+w_2y_l)]_q t}(-1)^{\lr (i_l+w_2 y_l)}.
\end{split}
\ee

From (\ref{10}), we note that

\be\la{11}
\begin{split}
&\sum_{\js}^{w_1-1}(-1)^{\lr j_l}q^{w_2 \lr(h-l)j_l} \I q^{w_1\lr (h-l)y_l}\\
&\hspace{2cm}\times e^{[w_1]_q[w_2x+\fa{w_2}{w_1}\lr j_l+\lr y_l]_{q^{w_1}}t}\dm\\
&=\LN \sum_{\js}^{w_1-1}\sum_{\is}^{w_2-1}\sum_{\ys}^{p^N-1}(-1)^{\lr(i_l+j_l+y_l)} q^{\lr(h-l)(w_2j_l+w_1i_l+w_1w_2y_l)}\\
&\hspace{2cm}\times e^{[w_1w_2(x+\lr y_l)+\lr w_2j_l+w_1 i_l]_q t}.
\end{split}
\ee

By the same method as (\ref{11}), we get

\be\la{12}
\begin{split}
&\sum_{\js}^{w_2-1}(-1)^{\lr j_l}q^{w_1 \lr(h-l)j_l} \I q^{w_2\lr (h-l)y_l}\\
&\hspace{2cm}\times e^{[w_2]_q[w_1x+\fa{w_1}{w_2}\lr j_l+\lr y_l]_{q^{w_2}}t}\dm\\
&=\LN \sum_{\js}^{w_2-1}\sum_{\is}^{w_1-1}\sum_{\ys}^{p^N-1}(-1)^{\lr(i_l+j_l+y_l)} q^{\lr(h-l)(w_1j_l+w_2i_l+w_1w_2y_l)}\\
&\hspace{2cm}\times e^{[w_1w_2(x+\lr y_l)+\lr w_1j_l+w_2 i_l]_q t}
\end{split}
\ee

Therefore, by (\ref{11}) and (\ref{12}), we obtain the following theorem.

\bt\la{t1}
For $w_1, w_2\in \N$ with $w_1\equiv 1\pmod{2}$ and $w_2\equiv 1\pmod{2}$, we have
\ben
\begin{split}
&\sum_{\js}^{w_1-1}(-1)^{\lr j_l}q^{w_2 \lr(h-l)j_l} \I q^{w_1\lr (h-l)y_l}\\
&\hspace{2cm}\times e^{[w_1]_q[w_2x+\fa{w_2}{w_1}\lr j_l+\lr y_l]_{q^{w_1}}t}\dm\\
&=\sum_{\js}^{w_2-1}(-1)^{\lr j_l}q^{w_1 \lr(h-l)j_l} \I q^{w_2\lr (h-l)y_l}\\
&\hspace{2cm}\times e^{[w_2]_q[w_1x+\fa{w_1}{w_2}\lr j_l+\lr y_l]_{q^{w_2}}t}\dm
\end{split}
\een
\et

\bc\la{c2}
For $n\geq 0$, $w_1, w_2\in \N$ with $w_1\equiv 1\pmod{2}$ and $w_2\equiv 1\pmod{2}$, we have
\ben
\begin{split}
&[w_1]_q^n\sum_{\js}^{w_1-1}(-1)^{\lr j_l}q^{w_2 \lr(h-l)j_l} \I q^{w_1\lr (h-l)y_l}\\
&\hspace{1cm}\times[w_2x+\fa{w_2}{w_1}\lr j_l+\lr y_l]_{q^{w_1}}^n \dm\\
&=[w_2]_q^n\sum_{\js}^{w_2-1}(-1)^{\lr j_l}q^{w_1 \lr(h-l)j_l} \I q^{w_2\lr (h-l)y_l}\\
&\hspace{1cm}\times[w_1x+\fa{w_1}{w_2}\lr j_l+\lr y_l]_{q^{w_2}}^n \dm\\
\end{split}
\een
\ec

Therefore, by (\ref{9}) and Corollary \ref{c2}, we obtain the following theorem.

\bt\la{t3}
For $n\geq 0$, $w_1, w_2\in \N$ with $w_1\equiv 1\pmod{2}$ and $w_2\equiv 1\pmod{2}$, we have
\ben
\begin{split}
&[w_1]_q^n\sum_{\js}^{w_1-1}(-1)^{\lr j_l}q^{w_2 \lr(h-l)j_l} E_{n,q^{w_1}}^{(h,r)}(w_2x+\fa{w_2}{w_1}(\JS))\\
&=[w_2]_q^n\sum_{\js}^{w_2-1}(-1)^{\lr j_l}q^{w_1 \lr(h-l)j_l}E_{n,q^{w_2}}^{(h,r)}(w_1x+\fa{w_1}{w_2}(\JS)). \\
\end{split}
\een
\et

From (\ref{3}), we note that
\be\la{13}
\begin{split}
&\I q^{w_1\lr(h-l)y_l}[w_2x+\fa{w_2}{w_1}\lr j_l+\lr y_l]_{q^{w_1}}^n\dm\\
&=\sum_{i=0}^{n}\bi{n}{i}\lp \fa{[w_2]_q}{[w_1]_q} \rp^i[\JS]_{q^{w_2}}^i q^{w_2(n-i)\lr j_l}\I [w_2x+\lr y_l]_{q^{w_1}}^{n-i}\\
&\hspace{4cm}\times q^{w_1\lr(h-l)y_l}\dm\\
&=\sum_{i=0}^{n}\bi{n}{i}\lp \fa{[w_2]_q}{[w_1]_q} \rp^i[\JS]_{q^{w_2}}^i q^{w_2(n-i)\lr j_l}E_{n-i,q^{w_1}}^{(h,r)}(w_2x).
\end{split}
\ee

By (\ref{13}), we get
\be\la{14}
\begin{split}
&[w_1]_q^n\sum_{\js}^{w_1-1}(-1)^{\lr j_l}q^{w_2 \lr(h-l)j_l} \I q^{w_1\lr (h-l)y_l}\\
&\hspace{4cm}\times [w_2x+\fa{w_2}{w_1}\lr j_l +\lr y_l]_{q^{w_1}}^n\dm\\
&=\sum_{\js}^{w_1-1}(-1)^{\lr j_l}q^{w_2 \lr(h-l)j_l}\sum_{i=0}^{n}\bi{n}{i}[w_2]_q^{i}[w_1]_q^{n-i}[\JS]_{q^{w_2}}^i \\
&\hspace{4cm}\times q^{w_2(n-1)\lr j_l}E_{n-i,q^{w_1}}^{(h,r)}(w_2x)\\
&=\sum_{i=0}^{n}\bi{n}{i}[w_2]_q^{i}[w_1]_q^{n-i}E_{n-i,q^{w_1}}^{(h,r)}(w_2x)\sum_{\js}^{w_1-1}(-1)^{\lr j_l}q^{w_2 \lr(n+h-l-i)j_l}\\
&\hspace{4cm}\times [\JS]_{q^{w_2}}^i\\
&=\sum_{i=0}^{n}\bi{n}{i}[w_2]_q^{i}[w_1]_q^{n-i}E_{n-i,q^{w_1}}^{(h,r)}(w_2x)T_{n,i,q^{w_2}}^{(h,r)}(w_1),
\end{split}
\ee
where
\be\la{15}
T_{n,i,q}^{(h,r)}(w)=\sum_{\js}^{w-1}(-1)^{\lr j_l}q^{ \lr(n+h-l-i)j_l}[\JS]_{q}^i.
\ee

By the same method as (\ref{14}), we get

\be\la{16}
\begin{split}
&[w_2]_q^n\sum_{\js}^{w_2-1}(-1)^{\lr j_l}q^{w_1 \lr(h-l)j_l} \I q^{w_2\lr (h-l)y_l}\\
&\hspace{2cm}\times[w_1x+\fa{w_1}{w_2}\lr j_l+\lr y_l]_{q^{w_2}}^n\dm\\
&=\sum_{\js}^{w_2-1}(-1)^{\lr j_l}q^{w_1 \lr(h-l)j_l}\sum_{i=0}^{n}\bi{n}{i}[w_1]_q^{i}[w_2]_q^{n-i}[\JS]_{q^{w_1}}^i\\
&\hspace{2cm}\times q^{w_1(n-i)\lr j_l}E_{n-i,q^{w_2}}^{(h,r)}(w_1x)\\
&=\sum_{i=0}^{n}\bi{n}{i}[w_1]_q^{i}[w_2]_q^{n-i}E_{n-i,q^{w_2}}^{(h,r)}(w_1x)T_{n,i,q^{w_1}}^{(h,r)}(w_2).
\end{split}
\ee

Therefore, by (\ref{14}) and (\ref{15}), we obtain the following theorem.

\bt\la{t4}
For $n\geq 0$, and  $w_1, w_2\in \N$ with $w_1\equiv 1\pmod{2}$ and\\
$w_2\equiv 1\pmod{2}$, we have
\ben
\begin{split}
&\sum_{i=0}^{n}\bi{n}{i}[w_2]_q^{i}[w_1]_q^{n-i}E_{n-i,q^{w_1}}^{(h,r)}(w_2x)T_{n,i,q^{w_2}}^{(h,r)}(w_1)\\
&=\sum_{i=0}^{n}\bi{n}{i}[w_1]_q^{i}[w_2]_q^{n-i}E_{n-i,q^{w_2}}^{(h,r)}(w_1x)T_{n,i,q^{w_1}}^{(h,r)}(w_2),
\end{split}
\een
where
\ben
T_{n,i,q}^{(h,r)}(w)=\sum_{\js}^{w-1}(-1)^{\lr j_l}q^{ \lr(n+h-l-i)j_l}[\JS]_{q}^i.
\een
\et



\begin{thebibliography}{99}

\bibitem {01} G. Andrews, B. C. Berndt, {\it Ramanujan's lost notebook, Part IV,} Springer, New York, 2013, xviii+439 pp.


\bibitem {02} S. Araci, M. Acikgoz, {\it A note on the Frobenius-Euler numbers and polynomials associated with Bernstein polynomials}, Adv. Stud. Contemp. Math., 22(2012),  no.3,  399-406.

\bibitem {03} I. N. Cangul, V. Kurt, H. Ozden, Y. Simsek, {\it On the higher-order $w$-$q$-Genocchi numbers}, Adv. Stud. Contemp. Math.  19(2009),  no. 1, 39-57.

\bibitem {04} D. S. Kim, N. Lee, J. Na, K. H. Park, {\it Identities of symmetry for higher-order Euler polynomials in three variables (I)}, Adv. Stud. Contemp. Math.  22(2012),  no. 1, 51-74.

\bibitem {05} T. Kim, {\it Barnes-type multiple $q$-zeta functions and $q$-Euler polynomials}, J. Phys. A 43(2010), no. 25, 255201, 11pp.

\bibitem {06} T. Kim, {\it Symmetry $p$-adic invariant integral on $\ma{Z}_p$ for Bernoulli and Euler polynomials}, J. Difference Equ. Appl. 14(2008), no.12, 1267-1277.

\bibitem {07} T. Kim, {\it Symmetry of power sum polynomials and multivariate fermionic $p$-adic invariant integral on $\zp$}, Russ. J. Math. Phys. 16(2009), no. 1, 93-96.

\bibitem {08} T. Kim, {\it An identity of the symmetry for the Frobenius-Euelr polynomials associated with the fermionic $p$-adic invariant $q$-integrals on $\zp$}, Rocky Mountain J. Math.  41(2011), no. 1, 239-247.

\bibitem {09} Y.-H. Kim, K.-W. Hwang, {\it Symmetry of power sum and twisted Bernoulli polynomials}, Adv. Stud. Contemp. Math. 18(2009),  no.2,  127-133.

\bibitem {10} H. Ozden, Y. Simsek, S.-H. Rim, I. N. Cangul, {\it A note on $p$-adic $q$-Euler measure}, Adv. Stud. Contemp. Math. 14(2007),  no.2,  233-239.

\bibitem {11} S. H. Rim, J. Jeong, {\it On the modified $q$-Euler numbers of higher order with weight}, Adv. Stud. Contemp. Math. 22(2012), 93-98.

\bibitem {12} Y. Simsek, {\it Interpolation functions of the Eulerian type polynomials and numbers}, Adv. Stud. Contemp. Math., 23(2013),  no.2, 301-307.

\bibitem {13}  Y. Simsek,  {\it Identities assoiciated with generalized Stirling type numbers and Eulerian type polynomials}, Math. Comput. Appl. 18(2013), no. 3, 251-263.

\bibitem {14}  H. J. H. Tuenter,  {\it A symmetry of power sum polynomials and Bernoulli numbers}, Amer. Math. Monthly 108(2001), no. 3, 258-261.


\end{thebibliography}
\end{document}